\documentclass[12pt]{article}
\usepackage{amsthm,amsfonts,amssymb,amscd}

\begin{document}


\begin{center}
\large \bf Birationally rigid Fano fibre spaces. II
\end{center}\vspace{0.5cm}

\centerline{A.V.Pukhlikov}\vspace{0.5cm}

\parshape=1
3cm 10cm \noindent {\small \quad\quad\quad \quad\quad\quad\quad
\quad\quad\quad {\bf }\newline In this paper we prove birational
rigidity of large classes of Fano-Mori fibre spaces over a base of
arbitrary dimension, bounded from above by a constant that depends
on the dimension of the fibre only. In order to do that, we first
show that if every fibre of a Fano-Mori fibre space satisfies
certain natural conditions, then every birational map onto another
Fano-Mori fibre space is fibre-wise. After that we construct large
classes of fibre spaces (into Fano double spaces of index one and
into Fano hypersurfaces of index one) which satisfy those
conditions.

Bibliography: 35 titles.} \vspace{1cm}

\section*{Introduction}

{\bf 0.1. Birationally rigid Fano-Mori fibre spaces.} In this
paper we investigate the problem of birational rigidity of
Fano-Mori fibre spaces $\pi\colon V\to S$. We assume that the base
$S$ is non-singular, the variety $V$ has at most factorial
terminal singularities, the anticanonical class $(-K_V)$ is
relatively ample and
$$
\mathop{\rm Pic} V = {\mathbb Z} K_V\oplus \pi^* \mathop{\rm Pic}
S.
$$
Let $\pi'\colon V'\to S'$ be an arbitrary rationally connected
fibre space, that is, a morphism of projective algebraic
varieties, where the base $S'$ and the fibre of general position
${\pi'}^{-1}(s')$, $s'\in S'$, are rationally connected and
$\mathop{\rm dim} V = \mathop{\rm dim} V'$. Consider a birational
map $\chi\colon V\dashrightarrow V'$ (provided they exist). In
order to describe the properties of the map $\chi$, of crucial
importance is whether $\chi$ is fibre-wise or not, that is,
whether this map transforms the fibres of the projection $\pi$
into the fibres of the projection $\pi'$. It is expected (and
confirmed by all known examples, see subsection 0.6), that the
answer is positive if the fibre space $\pi$ is ``sufficiently
twisted over the base''. Investigating this problem, one can
choose various classes of Fano-Mori fibre space and various
interpretations of the property to be ``twisted over the base''.
In the present paper we prove the following fact.\vspace{0.1cm}

{\bf Theorem 1.} {\it Assume that the Fano-Mori fibre space
$\pi\colon V\to S$ satisfies the conditions\vspace{0.1cm}

{\rm (i)} every fibre $F_s={\pi}^{-1}(s)$, $s\in S$, is a
factorial Fano variety with terminal singularities and the Picard
group $\mathop{\rm Pic} F_s = {\mathbb Z} K_{F_s}$,\vspace{0.1cm}

{\rm (ii)} for every effective divisor $D\in |-nK_{F_s}|$ on an
arbitrary fibre $F_s$ the pair $(F_s,\frac{1}{n} D)$ is log
canonical, and for every mobile linear system $\Sigma_s\subset
|-nK_{F_s}|$ the pair $(F_s,\frac{1}{n} D)$ is canonical for a
general divisor $D\in\Sigma_s$,\vspace{0.1cm}

{\rm (iii)} for every mobile family $\overline{\cal C}$ of curves
on the base $S$, sweeping out $S$, and a curve $\overline{C}\in
\overline{\cal C}$ the class of the following algebraic cycle of
dimension $\mathop{\rm dim} F$ for any positive $N\geqslant 1$
$$
-N (K_V\circ \pi^{-1}(\overline{C}))-F
$$
(where $F$ is the fibre of the projection $\pi$) is not effective,
that is, it is not rationally equivalent to an effective cycle of
dimension $\mathop{\rm dim} F$.\vspace{0.1cm}

Then every birational map $\chi\colon V\dashrightarrow V'$ onto
the total space of a rationally connected fibre space $V'/S'$ is
fibre-wise, that is, there exists a rational dominant map
$\beta\colon S\dashrightarrow S'$, such that the following diagram
commutes}
$$
\begin{array}{rcccl}
   & V & \stackrel{\chi}{\dashrightarrow} & V' & \\
\pi & \downarrow &   &   \downarrow & \pi' \\
   & S & \stackrel{\beta}{\dashrightarrow} & S'.
\end{array}
$$

Now we list the standard implications of Theorem 1, after which we
discuss the point of how restrictive the conditions (i)-(iii)
are.\vspace{0.1cm}

{\bf Corollary 1.} {\it In the assumptions of Theorem 1 on the
variety $V$ there are no structures of a rationally connected
fibre space over a base of dimension higher than $\mathop{\rm dim}
S$. In particular, the variety $V$ is non-rational. Every
birational self-map of the variety $V$ is fibre-wise and induces a
birational self-map of the base $S$, so that there is a natural
homomorphism of groups $\rho\colon \mathop{\rm Bir} V\to
\mathop{\rm Bir} S$, the kernel of which $\mathop{\rm Ker} \rho$
is the group $\mathop{\rm Bir} F_{\eta} = \mathop{\rm Bir} (V/S)$
of birational self-maps of the generic fibre $F_{\eta}$ (over the
generic non-closed point $\eta$ of the base $S$), whereas the
group $\mathop{\rm Bir} V$ is an extension of the normal subgroup
$\mathop{\rm Bir} F_{\eta}$ by the group $\Gamma=\rho(\mathop{\rm
Bir} V)\subset \mathop{\rm Bir} S$:}
$$
1\to \mathop{\rm Bir} F_{\eta} \to \mathop{\rm Bir} V\to \Gamma
\to 1.
$$

How restrictive are the conditions (i)-(iii)? The condition (iii)
belongs to the same class of conditions as the well known
$K^2$-condition and the $K$-condition for fibrations over
${\mathbb P}^1$ (see, for instance, \cite[Chapter 4]{Pukh13a}) and
the Sarkisov condition for conic bundles (see \cite{S80,S82}).
This condition measures the ``degree of twistedness'' of the fibre
space $V/S$ over the base $S$. Below we illustrate this meaning of
the condition (iii) by particular examples. We will see that this
condition is not too restrictive: for a fixed method of
constructing the fibre space $V/S$ and a fixed ``ambient'' fibre
space $X/S$ the condition (iii) is satisfied by ``almost all''
families of fibre spaces $V/S$.\vspace{0.1cm}

In terms of numerical geometry of the varieties $V$ and $S$ the
condition (iii) can be expressed in the following way. Let
$$
A^*(V)=\mathop{\bigoplus}\limits^{\mathop{\rm dim} V}_{i=0} A^i(V)
$$
be the numerical Chow ring of the variety $V$, graded by
codimension. Set
$$
A_i(V)=A^{\mathop{\rm dim} V -i} (V) \otimes {\mathbb R}
$$
and denote by the symbol $A^{\rm mov}_i (V)$ the closed cone in
$A_i(V)$, generated by the classes of mobile cycles, the families
of which sweep out $V$, and by the symbol $A^+_i (V)$ the
pseudoeffective cone in $A_i(V)$, generated by the classes of
effective cycles. Furthermore, by the symbol $A_{i,\leqslant
j}(V)$ we denote the linear subspace in $A_i(V)$, generated by the
classes of subvarieties of dimension $i$, the image of which on
$S$ has dimension at most $j$. In the real space $A_{i,\leqslant
j}(V)$ consider the closed cones $A_{i,\leqslant j}^{\rm mov}(V)$
of mobile and $A_{i,\leqslant j}^+(V)$ of pseudoeffective classes.
In a similar way we define the real vector space $A_i(S)$ and the
closed cones $A^{\rm mov}_i (S)$ è $A_i^+(S)$. If
$\delta=\mathop{\rm dim} F$ is the dimension of the fibre of the
projection $\pi$, then the operation of taking the preimage
generates a linear map
$$
\pi^*A_i(S) \to A_{\delta+i,\leqslant i}(V),
$$
whereas $\pi^*(A_i^+(S))\subset A_{\delta+i,\leqslant i}^+(V)$ and
$\pi^*(A_i^{\rm mov}(S))\subset A_{\delta+i,\leqslant i}^{\rm
mov}(V)$. Now let us consider the linear map
$$
\gamma\colon A_1(S)\to A_{\delta,\leqslant 1}(V),
$$
defined by the formula
$$
z\mapsto -(K_V\cdot \pi^* z).
$$
The condition (iii) means that the image of the cone
$\gamma(A_1^{\rm mov} (S))$ is contained in the boundary of the
pseudoeffective cone $A^+_{\delta,\leqslant 1}(V)$, that is,
$$
\gamma (A_1^{\rm mov} (S)) \cap \mathop{\rm Int}
A^+_{\delta,\leqslant 1}(V) = \emptyset.
$$
More precisely, for any class $z\in A_1^{\rm mov} (S)$ the
intersection of the closed ray
$$
\{ \gamma(z) - t[F]\,|\, t\in {\mathbb R}_+ \}
$$
(where $[F]\in \mathop{\rm Int} A^+_{\delta,\leqslant 1}(V)$ is
the class of the fibre of the projection $\pi$) with the cone
$A^+_{\delta,\leqslant 1}(V)$ either is empty or consists of just
one point $\gamma(z)$.\vspace{0.1cm}

One may suggest that the condition (iii) is close to a precise one
(``if and only if''), that is, its violation (or an essential
deviation from this condition) implies the existence of another
structure of a Fano-Mori fibre space on the variety
$V$.\vspace{0.1cm}

The following remark gives an obvious way to check the condition
(iii).\vspace{0.1cm}

{\bf Remark 0.1.} Assume that on the variety $V$ there is a
numerically effective divisorial class $L$ such that
$(L^{\delta}\cdot F)>0$ and the linear function $(\cdot
L^{\delta})$ is non-positive on the cone $\gamma(A_1^{\rm mov}
(S))$, that is to say, for any mobile curve $\overline{C}$ on $S$
the inequality
\begin{equation}\label{10.06.2014.1}
\left(L^{\delta}\cdot K_V\cdot
\pi^{-1}(\overline{C})\right)\geqslant 0
\end{equation}
holds. Then the condition (iii) is obviously
satisfied.\vspace{0.1cm}

The conditions (i) and (ii), however, are much more restrictive.
They mean that {\it all} fibres of the projection $\pi$ are
varieties of sufficiently general position in their family. This
implies that the dimension of the base for a fixed family of
fibres is bounded from above (by a constant depending on the
particular family, to which the fibres belong). In the examples
considered in the present paper, for a sufficiently high dimension
of the fibre $\delta=\mathop{\rm dim} F$ the dimension of the base
is bounded from above by a number of order $\frac12\delta^2$.
Recall that up to now not a single example was known of a
fibration into higher dimensional Fano varieties over a base of
dimension two and higher with just one structure of a rationally
connected fibre space (for a brief historical survey, see
subsection 0.5).\vspace{0.3cm}


{\bf 0.2. Fibrations into double spaces of index one.} By the
symbol ${\mathbb P}$ we denote the projective space ${\mathbb
P}^M$, $M\geqslant 5$. Let ${\cal W}={\mathbb P}(H^0({\mathbb
P},{\cal O}_{\mathbb P}(2M)))$ be the space of hypersurfaces of
degree $2M$ in ${\mathbb P}$. The following fact is
true.\vspace{0.1cm}

{\bf Theorem 2.} {\it There exists a Zariski open subset ${\cal
W}_{\rm reg}\subset {\cal W}$, such that for any hypersurface
$W\in {\cal W}_{\rm reg}$ the double cover $\sigma\colon F\to
{\mathbb P}$, branched over $W$, satisfies the conditions (i) and
(ii) of Theorem 1, and moreover, the estimate}
$$
\mathop{\rm codim} (({\cal W}\setminus {\cal W}_{\rm reg})\subset
{\cal W})\geqslant \frac{(M-4)(M-1)}{2}
$$
{\it holds.}\vspace{0.1cm}

An explicit description of the set ${\cal W}_{\rm reg}$ and a
proof of Theorem 2 are given in \S 2. Fix a number $M\geqslant 5$
and a non-singular rationally connected variety $S$ of dimension
$\mathop{\rm dim} S<\frac12 (M-4)(M-1)$. Let ${\cal L}$ be a
locally free sheaf of rank $M+1$ on $S$ and $X={\mathbb P}({\cal
L})={\bf Proj}\, \mathop{\oplus}\limits_{i=0}^{\infty}{\cal
L}^{\otimes i}$ the corresponding ${\mathbb P}^M$-bundle. We may
assume that ${\cal L}$ is generated by its sections, so that the
sheaf ${\cal O}_{{\mathbb P}({\cal L})}(1)$ is also generated by
the sections. Let $L\in \mathop{\rm Pic} X$ be the class of that
sheaf, so that
$$
\mathop{\rm Pic} X = {\mathbb Z} L\oplus \pi_X^* \mathop{\rm Pic}
S,
$$
where $\pi_X\colon X\to S$ is the natural projection. Take a
general divisor $U\in |2(ML+\pi^*_X R)|$, where $R\in \mathop{\rm
Pic} S$ is some class. If this system is sufficiently mobile, then
by the assumption about the dimension of the base $S$ and Theorem
2 we may assume that for every point $s\in S$ the hypersurface
$U_s= U\cap \pi^{-1}_X(s)\in {\cal W}_{\rm reg}$, and for that
reason the double space branched over $U_s$, satisfies the
conditions (i) and (ii) of Theorem 1. Let $\sigma\colon V\to X$
the double cover branched over $U$. Set
$\pi=\pi_X\circ\sigma\colon V\to S$, so that $V$ is a fibration
into Fano double spaces of index one over $S$. Recall that the
divisor $U\in |2(ML+\pi^*_X R)|$ is assumed to be sufficiently
general.\vspace{0.1cm}

{\bf Theorem 3.} {\it Assume that the divisorial class $(K_S+R)$
is pseudoeffective. Then for the fibre space $\pi\colon V\to S$
the claims of Theorem 1 and Corollary 1 hold. In particular,
$$
\mathop{\rm Bir} V = \mathop{\rm Aut} V = {\mathbb Z}/2{\mathbb Z}
$$
is the cyclic group of order 2.}\vspace{0.1cm}

{\bf Proof.} Since the conditions (i) and (ii) of Theorem 1 are
satisfied by construction of the variety $V$, it remains to check
the condition (iii). Let us use Remark 0.1. Elementary
computations show that the inequality (\ref{10.06.2014.1}) up to a
positive factor is the inequality
$$
((K_S+R)\cdot \overline{C})\geqslant 0.
$$
Since the curve $\overline{C}$ belongs to a mobile family,
sweeping out the base $S$, the last inequality holds if the class
$(K_S+R)$ is pseudoeffective. Q.E.D. for the
theorem.\vspace{0.1cm}

{\bf Example 0.2.} Take $S = {\mathbb P}^m$, where $m<\frac12
(M-4)(M-1)$, $X={\mathbb P}^M\times {\mathbb P}^m$ and $W_X$ is a
generic hypersurface of bidegree $(2M,2l)$, where $l\geqslant
m+1$. Then for the double cover $\sigma\colon V\to X$, branched
over $W_X$, the claims of Theorem 1 and Corollary 1 are true. Note
that for $l\leqslant m$ on the double cover $V$ there is another
structure of a Fano fibre space: it is given by the projection
$\pi_1\colon V\to {\mathbb P}^M$. Therefore, the condition (iii)
of Theorem 1 and its realization in Theorem 3 turn out to be
precise.\vspace{0.3cm}


{\bf 0.3. Fibrations into Fano hypersurfaces of index one.} The
symbol ${\mathbb P}$ still stands for the projective space
${\mathbb P}^M$, $M\geqslant 10$. Fix $M$. Let ${\cal F}={\mathbb
P}(H^0({\mathbb P},{\cal O}_{\mathbb P}(M)))$ be the space of
hypersurfaces of degree $M$ in ${\mathbb P}$. The following fact
is true.\vspace{0.1cm}

{\bf Theorem 4.} {\it There is a Zariski open subset ${\cal
F}_{\rm reg}\subset {\cal F}$, such that every hypersurface $F\in
{\cal F}_{\rm reg}$ satisfies the conditions (i) and (ii) of
Theorem 1, and the following estimate holds:}
\begin{equation}\label{14.06.2014.1}
\mathop{\rm codim} (({\cal F}\setminus {\cal F}_{\rm reg})\subset
{\cal F})\geqslant \frac{(M-7)(M-6)}{2}-5.
\end{equation}

An explicit description of the subset ${\cal F}_{\rm reg}$ and a
proof of Theorem 4 are given in \S 2-3. Fix a non-singular
rationally connected variety $S$ of dimension $\mathop{\rm dim}
S<\frac12 (M-7)(M-6)-5$. As in subsection 0.2, let ${\cal L}$ be a
locally free sheaf of rank $M+1$ on $S$ and $X={\mathbb P}({\cal
L})={\bf Proj}\, \mathop{\oplus}\limits_{i=0}^{\infty}{\cal
L}^{\otimes i}$ the corresponding ${\mathbb P}^M$-bundle in the
sense of Grothendieck; we assume that ${\cal L}$ is generated by
global sections. Let $\pi_X\colon X\to S$ be the projection, $L\in
\mathop{\rm Pic} X$ the class of the sheaf ${\cal O}_{{\mathbb
P}({\cal L})}(1)$. Consider a general divisor
$$
V\in |ML+\pi^*_X R|,
$$
where $R\in \mathop{\rm Pic} S$ is some divisor on the base. By
the assumption about the dimension of the base made above and
Theorem 4 we may assume that the Fano fibre space $\pi\colon V\to
S$, where $\pi=\pi_X|_V$, satisfies the conditions (i) and (ii) of
Theorem 1.\vspace{0.1cm}

{\bf Theorem 5.} {\it Assume that the divisorial class
$(K_S+\left(1-\frac{1}{M}\right)R)$ is pseudoeffective. Then for
the Fano fibre space $\pi\colon V\to S$ the claims of Theorem 1
and Corollary are true. In particular, the group
$$
\mathop{\rm Bir} V = \mathop{\rm Aut} V
$$
is trivial.}\vspace{0.1cm}

{\bf Proof.} The conditions (i) and (ii) of Theorem 1 are
satisfied by the generality of the divisor $V$. The inequality
(\ref{10.06.2014.1}) up to a positive factor is the same as the
inequality
$$
((M K_S+(M-1)R)\cdot \overline{C})\geqslant 0.
$$
Therefore, by Remark 0.1, the condition (iii) of Theorem 1 also
holds. Q.E.D. for the theorem.\vspace{0.1cm}

{\bf Example 0.2.} Take $S = {\mathbb P}^m$, where
$m\leqslant\frac12 (M-7)(M-6)-6$, $X={\mathbb P}^M\times {\mathbb
P}^m$ and $V\subset X$ is a sufficiently general hypersurface of
bidegree $(M,l)$, where $l$ satisfies the inequality
$$
l\geqslant\frac{M}{M-1} (m+1).
$$
Then the Fano fibre space $V/{\mathbb P}^m$ satisfies all
assumptions of Theorem 1 and therefore for this fibre space the
claim of Theorem 1 and that of Corollary 1 are true. Note that for
$l\leqslant m$ on the variety $V$ there is another structure of a
Fano fibre space, given by the projection $V\to {\mathbb P}^M$.
Note also that if we fix the dimension $m$ of the base, then for
$M\geqslant m$ the condition of Theorem 5 is close to the optimal
one: it is satisfied for $l\geqslant m+2$, so that the only value
of the integral parameter $l$, for which the problem of birational
rigidity of the fibre space $V/{\mathbb P}^m$ remains open, is
$l=m+1$. In that case the projection $V\to {\mathbb P}^M$ is a
$K$-trivial fibre space.\vspace{0.3cm}


{\bf 0.4. The structure of the paper.} The present paper is
organized in the following way. In \S 1 we prove Theorem 1. After
that, in \S 2 we deal with the conditions of general position,
which should be satisfied for every fibre of the fibre space $V/S$
in order for the conditions (i) and (ii) of Theorem 1 to hold. The
conditions of general position (regularity) are given for Fano
double spaces of index one and Fano hypersurfaces of index one.
This makes it possible to define the sets ${\cal W}_{\rm reg}$ and
${\cal F}_{\rm reg}$ and prove Theorem 2 and carry out the
preparational work for the proof of Theorem 4, the main technical
fact of the present paper, which implies Theorem 5, geometrically
the most impressive result of this paper, in an obvious
way.\vspace{0.1cm}

In \S 3 we complete the proof of Theorem 4, more precisely we show
that the condition (ii) of Theorem 1 is satisfied for a regular
Fano hypersurface $F\in {\cal F}_{\rm reg}$. The proof makes use a
combination of the technique of hypertangent divisors and the
inversion of adjunction. Note that the approach of the present
paper corresponds to the linear method of proving birational
rigidity, see \cite[Chapter 7]{Pukh13a}; the technique of the
quadratic method (in the first place, the technique of counting
multiplicities) is not used.\vspace{0.1cm}

The assumption in Theorem 1 that the base $S$ of the fibre space
$V/S$ is non-singular seems to be unnecessary and could be
replaced by the condition that the singularities are at most
terminal and (${\mathbb Q}$-)factorial.\vspace{0.3cm}


{\bf 0.5. Historical remarks and acknowledgements.} The starting
point of studying birational geometry of rationally connected
fibre spaces seems to be the use of de Jonqiere transformations
(see, for instance, \cite{Hud}). In the modern algebraic geometry
the objects of this type started to be systematically investigated
in the works of V.A.Iskovskikh and M.Kh.Gizatullin about pencils
of rational curves \cite{I67,I70,Giz67} over non-closed fields,
which followed the investigation of the ``absolute'' case in the
papers of Yu.I.Manin \cite{M66,M67,M72}. We also point out the
paper of I.V.Dolgachev \cite{Dol66}, which started (in the modern
period) the study of $K$-trivial fibrations.\vspace{0.1cm}

After the breakthrough in three-dimensional birational geometry
that was made in the classical paper of V.A.Iskovskikh and
Yu.I.Manin on the three-dimensional quartic \cite{IM} the problems
of the ``relative'' three-dimensional birational geometry were the
next to be investigated, that is, the task was to describe
birational maps of three-dimensional algebraic varieties, fibred
into conics over a rational surface or into del Pezzo surfaces
over ${\mathbb P}^1$. The famous Sarkisov theorem gave an almost
complete solution of the question of birational rigidity for conic
bundles \cite{S80,S82}. A similar question for the pencils of del
Pezzo surfaces remained absolutely open until 1996 \cite{Pukh98a};
see the introduction to the last paper about the reasons of those
difficulties (the test class construction turned out to be
unsuitable for studying the varieties of that type).\vspace{0.1cm}

The method of proving birational rigidity, realized in
\cite{Pukh98a}, generalized well into the arbitrary dimension, for
varieties fibred into Fano varieties over ${\mathbb P}^1$. In a
long series of papers
\cite{Pukh00d,Sob01,Sob02,Pukh04a,Pukh04b,Pukh06a,Pukh06b,Pukh09a}
birational rigidity was shown for many classes of Fano fibre space
over ${\mathbb P}^1$. At the same time, the birational geometry of
the remaining families of three-dimensional varieties with a
pencil of del Pezzo surfaces of degree 1 and 2 was investigated
\cite{Grin00,Grin03a,Grin03b,Grin04}; in that direction the
results that were obtained were nearly exhaustive. However, the
base of the fibre spaces under investigation remained
one-dimensional and even Fano fibrations over surfaces seemed to
be out of reach.\vspace{0.1cm}

The only exception in that series of results was the theorem about
Fano direct products \cite{Pukh05} and the papers about direct
products that followed \cite{Pukh08a,Ch08}. In those papers the
Fano fibre spaces under consideration had the both the base and
the fibre of arbitrary dimension. However, the fibre spaces
themselves were very special (direct products) and could not
pretend to be  {\it typical} Fano fibre spaces.\vspace{0.1cm}

The present paper gives, at long last, numerous examples of
typical birationally rigid Fano fibre spaces with the base and
fibre of high dimension (for a fixed dimension of the fibre
$\delta$ the dimension of the base is bounded by a constant of
order $\sim \frac12\delta^2$). Theorem 1 can be viewed as a
realization of the well known principle: the ``sufficient
twistedness'' of a fibre space over the base implies birational
rigidity. This principle was many times confirmed in the class of
fibrations over ${\mathbb P}^1$; now it is extended to the the
class of fibre spaces over a base of arbitrary
dimension.\vspace{0.1cm}

The main object of study in this paper is a fibre space into Fano
hypersurfaces of index one, so that it is a follow up of the paper
\cite{Pukh00d}. From the technical viewpoint, the predecessors of
this paper are \cite{Pukh08a,Pukh09b}, where the {\it linear}
method of proving birational rigidity was developed. It is
possible, however, that the quadratic techniques could be applied
to the class of Fano fibre spaces over a base of arbitrary
dimension as well.\vspace{0.1cm}

Various technical moments related to the arguments of the present
paper were discussed by the author in his talks given in 2009-2014
at Steklov Institute of Mathematics. The author is grateful to the
members of the divisions of Algebraic Geometry and Algebra and
Number Theory for the interest in his work. The author also thanks
his colleagues in the Algebraic Geometry research group at the
University of Liverpool for the creative atmosphere and general
support.


\section{Birationally rigid fibre spaces}

In this section we prove Theorem 1. We do it in three steps:
first, assuming that the birational map $\chi\colon
V\dashrightarrow V'$ is not fibre-wise, we prove the existence of
a maximal singularity of the map $\chi$, covering the base $S'$
(subsection 1.1). After that, we construct such a sequence of blow
ups of the base $S^+\to S$, that the image of every maximal
singularity on $S$ is a prime divisor (subsection 1.2). Finally,
using a very mobile family of curves contracted by the projection
$\pi'$, we obtain a contradiction with the condition (iii) of
Theorem 1 (subsection 1.3). This implies that the map $\chi\colon
V\dashrightarrow V'$ is fibre-wise, which completes the proof of
Theorem 1.\vspace{0.3cm}

{\bf 1.1. Maximal singularities of birational maps.} In the
notations of Theorem 1 fix a birational map $\chi\colon
V\dashrightarrow V'$ onto the total space $V'$ of a rationally
connected fibre space $\pi'\colon V'\to S'$. Consider any very
ample linear system $\overline{\Sigma'}$ on $S'$. Let
$\Sigma'=(\pi')^*\overline{\Sigma'}$ be its pull back on $V'$, so
that the divisors $D'\in\Sigma'$ are composed from the fibres of
the projection $\pi'$, and for that reason for any curve $C\subset
V'$, contracted by the projection $\pi'$ we have $(D'\cdot C)=0$.
The linear system $\Sigma'$ is obviously mobile. Let
$$
\Sigma=(\chi^{-1})_*\Sigma'\subset |-nK_V+\pi^*Y|
$$
be its strict transform on $V$, where $n\in{\mathbb Z}_+$.
\vspace{0.1cm}

{\bf Lemma 1.1.} {\it For any mobile family of curves
$\overline{C}\in\overline{\cal C}$ on $S$, sweeping out $S$, the
inequality $(\overline{C}\cdot Y)\geqslant 0$ holds, that is to
say, the numerical class of the divisor $Y$ is non-negative on the
cone} $A^{\rm mov}_1(S)$.\vspace{0.1cm}

{\bf Proof.} This is almost obvious. For a general divisor
$D\in\Sigma$ the cycle $(D\circ \pi^{-1}(\overline{C}))$ is
effective. Its class is $-n(K_V\circ \pi^{-1}(\overline{C}))
+(Y\cdot \overline{C})F$, so that by the condition (iii) the claim
of the lemma follows. Q.E.D.\vspace{0.1cm}

Obviously, the map $\chi$ is fibre-wise if and only if $n=0$.
Therefore, if $n=0$, then the claim of Theorem 1 holds. So let us
assume that $n\geq 1$ and show that this assumption leads to a
contradiction.\vspace{0.1cm}

The linear system $\Sigma$ is mobile. Let us resolve the
singularities of the map $\chi$: let
$$
\varphi\colon\widetilde{V}\to V
$$
be a birational morphism (a composition of blow ups with
non-singular centres), where $\widetilde{V}$ is non-singular and
the composition
$\chi\circ\varphi\colon\widetilde{V}\dashrightarrow V'$ is
regular. Furthermore, consider the set ${\cal E}$ of prime
divisors on $\widetilde{V}$, satisfying the following conditions:

\begin{itemize}

\item every divisor $E\in{\cal E}$ is $\varphi$-exceptional,

\item for every $E\in{\cal E}$ the closed set $\chi\circ\varphi(E)
\subset V'$ is a prime divisor on $V'$,

\item the set $\chi\circ\varphi(E)$ for every
$E\in{\cal E}$ covers the base: $\pi'[\chi\circ\varphi(E)]=S'$.

\end{itemize}

Setting $\widetilde{K}=K_{\widetilde{V}}$, write down
$$
\widetilde{\Sigma}\subset|-n\widetilde{K}+(\pi^*Y- \sum_{E\in{\cal
E}}\varepsilon(E)E)+\Xi|,
$$
where $\widetilde{\Sigma}$, as usual, is the strict transform of
the mobile linear system $\Sigma$ on $\widetilde{V}$,
$\varepsilon(E)\in{\mathbb Z}$ is some coefficient and $\Xi$
stands for a linear combination of $\varphi$-exceptional divisors
which do not belong to the set ${\cal E}$.\vspace{0.1cm}

{\bf Definition 1.1.} An exceptional divisor $E\in{\cal E}$ is
called {\it a maximal singularity} of the map $\chi$, if
$\varepsilon(E)>0$.\vspace{0.1cm}

Obviously, a maximal singularity satisfies {\it the Noether-Fano
inequality}
$$
\mathop{\rm ord}\nolimits_E\varphi^*\Sigma>na(E),
$$
where $a(E)=a(E,V)$ is the discrepancy of the divisor $E$ with
respect to $V$. In this paper we somewhat modify the standard
concept of a maximal singularity, requiring in addition that it is
realized by a divisor on $V'$, covering  the base. Let ${\cal
M}\subset{\cal E}$ be the set of all maximal
singularities.\vspace{0.1cm}

{\bf Proposition 1.1.} {\it Maximal singularities do exist:}
${\cal M\neq\emptyset}$.\vspace{0.1cm}

{\bf Proof.} Assume the converse, that is, for any $E\in{\cal E}$
the inequality $\varepsilon(E)\leq 0$ holds. Let ${\cal C}'$ be a
family of rational curves on $V'$, satisfying the following
conditions:

\begin{itemize}

\item the curves $C'\in{\cal C}'$ are contracted by the projection $\pi'$,

\item the curves $C'\in{\cal C}'$ sweep out a dense open subset in $V'$,

\item the curves $C'\in{\cal C}'$ do not intersect the set of
points where the rational map $(\chi\circ\varphi)^{-1}\colon
V'\dashrightarrow\widetilde{V}$ is not well defined.

\end{itemize}

Apart from that, we assume that a general curve $C'\in{\cal C}'$
intersects every divisor $\chi\circ\varphi(E)$, $E\in{\cal E}$,
transversally at points of general position. Such a family of
curves we will call {\it very mobile}. Obviously, very mobile
families of rational curves do exist.\vspace{0.1cm}

Let $\widetilde{C}\cong C'$ be the inverse image of the curve
$C'\in{\cal C}'$ on $\widetilde{V}$. Since the linear system
${\Sigma}'$ is pulled back from the base, for a divisor
$\widetilde{D}\in\widetilde{\Sigma}$ we have the equality
$(\widetilde{C}\cdot\widetilde{D})=0$.  On the other hand,
$(\widetilde{C}\cdot\widetilde{K})=(C'\cdot K_{V'})<0$ and
$$
(\widetilde{C}\cdot(\pi^*Y-\sum_{E\in{\cal
E}}\varepsilon(E)E))\geqslant 0,
$$
since by the condition (iii) of our theorem $(\widetilde{C}
\cdot\pi^*Y)\geqslant 0$ and by assumption
$-\varepsilon(E)\in{\mathbb Z}_+$ for all $E\in{\cal E}$. Finally,
the divisor $\Xi$ (which is not necessarily effective) is a linear
combination of such $\varphi$-exceptional divisors
$R\subset\widetilde{V}$, that $\pi'[\chi\circ\varphi(R)]$ is a
proper closed subset of the base $S'$. So we have the equality
$(\widetilde{C}\cdot\Xi)=0$. This implies that
$$
(\widetilde{C}\cdot\widetilde{D})\geqslant n>0,
$$
which is a contradiction. Therefore, ${\cal M}\neq\emptyset$.
Q.E.D. for the proposition.\vspace{0.1cm}

{\bf Proposition 1.2.} {\it For any maximal singularity
$E\subset{\cal M}$ its center
$$
\mathop{\rm centre}(E,V)=\varphi(E)
$$
on $V$ does not cover the base: $\pi(\mathop{\rm
centre}(E,V))\subset S$ is a proper closed subset of the variety
$S$.}\vspace{0.1cm}

{\bf Proof.} Assume the converse: the centre of some maximal
singularity $E\in{\cal M}$ covers the base: $\pi(\mathop{\rm
centre}(E,V))=S$. Let $F=\pi^{-1}(s)$, $s\in S$ be a fibre of
general position. By assumption the strict transform
$\widetilde{F}$ of the fibre $F$ on $\widetilde{V}$ has a
non-empty intersection with $E$, and for that reason every
irreducible component of the intersection $\widetilde{F}\cap E$ is
a maximal singularity of the mobile linear system
$\Sigma_F=\Sigma|_F\subset |-nK_F|$. However, by the condition
(ii) of Theorem 1 on the variety $F$ there are no mobile linear
systems with a maximal singularity. This contradiction proves the
proposition.\vspace{0.3cm}


{\bf 1.2. The birational modification of the base of the fibre
space $V/S$.} Now let us construct a sequence of blow ups of the
base, the composition of which is a birational morphism
$\sigma_S\colon S^+\to S$, and the corresponding sequence of blow
ups of the variety $V$, the composition of which is a birational
morphism $\sigma\colon V^+\to V$, where
$V^+=V\mathop{\times}\nolimits_S S^+$, so that the following
diagram commutes
$$
\begin{array}{rcccl}
& V^+ &\stackrel{\sigma}{\to} & V &\\
\pi_+ & \downarrow & & \downarrow & \pi\\
& S^+ & \stackrel{\sigma_S}{\to} & S. &\\
\end{array}
$$
The birational morphism $\sigma_S$ is constructed inductively as a
composition of elementary blow ups $\overline{\sigma}_i\colon
S_i\to S_{i-1}$, $i=1,\dots$, where $S_0=S$. Assume that
$\overline{\sigma}_i$ are already constructed for $i\leqslant k$
(if $k=0$, then we start with the base $S$). Set
$V_k=V\mathop{\times}\nolimits_SS_k$ and let $\pi_k\colon V_k\to
S_k$ be the projection. Consider the irreducible closed subsets
\begin{equation}\label{20.05.2014.1}
\pi_k(\mathop{\rm centre}(E,V_k))\subset S_k,
\end{equation}
where $E$ runs through the set ${\cal M}$. By Proposition 1.2, all
these subsets are proper subsets of the base $S_k$. If all of them
are prime divisors on $S_k$, we stop the procedure: set $S^+=S_k$
and $V^+=V_k$. Otherwise, for $\overline{\sigma}_{k+1}$ we take
the blow up of any inclusion-minimal set (\ref{20.05.2014.1}) for
all $E\in{\cal M}$.\vspace{0.1cm}

It is easy to check that the sequence of blow ups
$\overline{\sigma}$ terminates. Indeed, set
$$
\alpha_k=\sum_{E\in{\cal M}}a(E,V_k).
$$
Since the birational morphism $\sigma_k\colon V_k\to V_{k-1}$ is
the blow up of a closed irreducible subset, containing the centre
of one of the divisors $E\in{\cal M}$ on $V_{k-1}$, we get the
inequality $\alpha_{k+1}<\alpha_k$. The numbers $\alpha_i$ are by
construction non-negative, which implies that the sequence of blow
ups $\overline{\sigma}_i$ is finite. Therefore, for any maximal
singularity $E\in{\cal M}$ the closed subset $\pi_+(\mathop{\rm
centre}(E,V^+))\subset S^+$ is a prime divisor.\vspace{0.3cm}


{\bf 1.3. The mobile family of curves.} Again let us consider a
very mobile family of curves ${\cal C}'$ on $V'$ and its strict
transform ${\cal C}^+$ on $V^+$. Let $C^+\in{\cal C}^+$ be a
general curve and $\overline{C^+}=\pi_+(C^+)$ the corresponding
curve of the family $\overline{\cal C^+}$ on $S^+$. Furthermore,
let $\Sigma^+$ be the strict transform of the linear system
$\Sigma$ on $V^+$. For some class of divisors $Y^+$ on $S^+$ we
have:
$$
\Sigma^+\subset |-nK^++\pi^*_+Y^+|,
$$
where for simplicity of notation $K^+=K_{V^+}$. Note that even if
$Y$ is an effective or mobile class on $S$, in this case $Y^+$ is
not its strict transform on $S^+$, that is to say, we violate the
principle of notations. The following observation is
crucial.\vspace{0.1cm}

{\bf Proposition 1.3.} {\it The inequality
$$
(\overline{C^+}\cdot Y^+)<0
$$
holds. In particular, the class $Y^+$ is not
pseudoeffective.}\vspace{0.1cm}

{\bf Proof.} Assume the converse:
$$
(C^+\cdot\pi^*Y^+)=(\overline{C^+}\cdot Y^+)\geqslant 0.
$$
We may assume that the resolution of singularities $\varphi$ of
the map $\chi$ filters through the sequence of blow ups
$\sigma\colon V^+\to V$, so that for the strict transform
$\widetilde{\Sigma}$ of the linear system $\Sigma$ on
$\widetilde{V}$ we have
$$
\widetilde{\Sigma}\subset |-n\widetilde{K}+(\pi^*_+Y^+
-\sum_{E\in{\cal E}}\widetilde{\varepsilon}(E)E)+\widetilde{\Xi}|,
$$
where $\widetilde{K}=K_{\widetilde{V}}$,
$\widetilde{\varepsilon}(E)\in{\mathbb Z}$ and $\widetilde{\Xi}$
is a linear combination of exceptional divisors of the birational
morphism $\widetilde{V}\to V^+$, which are not in the set ${\cal
E}$. For the strict transform $\widetilde{C}\in\widetilde{\cal C
}$ of the curve $C^+\in{\cal C}^+$ and the divisor
$\widetilde{D}\in\widetilde{\Sigma}$ we have, as in the proof of
Proposition 1.1, the equality
$(\widetilde{C}\cdot\widetilde{D})=0$. By the construction of the
divisor $\widetilde{\Xi}$ we have
$(\widetilde{C}\cdot\widetilde{\Xi})=0$. Finally,
$(\widetilde{C}\cdot\widetilde{K})<0$, whence we conclude that
$$
(\widetilde{C}\cdot(\pi^*_+Y^+-\sum_{E\in{\cal E}}
\widetilde{\varepsilon}(E)E)< 0.
$$

By our assumption for at least one divisor $E\in{\cal E}$ we have
the inequality $\widetilde{\varepsilon}(E)> 0$. This divisor is
automatically a maximal singularity, $E\in{\cal M}$. By our
construction, however, we can say more: $E$ is a maximal
singularity for the mobile linear system $\Sigma^+$ as well, that
is, the pair $\left(V^+,\frac{1}{n}\Sigma^+\right)$ is not
canonical and $E$ realizes a non-canonical singularity of that
pair.\vspace{0.1cm}

However, $\pi_+(\mathop{\rm centre}(E,V^+))=\overline{E}\subset
S^+$ is a prime divisor, so that $\pi^{-1}_+(\overline{E})\subset
V^+$ is also a prime divisor. The linear system $\Sigma^+$ has no
fixed components, therefore for a general point $s\in\overline{E}$
and the corresponding fibre $F=\pi^{-1}_+(s)\subset V^+$ we have:
the linear system $\Sigma_F=\Sigma^+|_F\subset |-nK_F|$ is
non-empty and for $D_F\in\Sigma_F$ the pair
$\left(F,\frac{1}{n}D_F\right)$ is non log canonical by the
inversion of adjunction (see \cite{Kol93}). This contradicts the
condition (ii) of our theorem. Proposition 1.3 is shown.
Q.E.D.\vspace{0.1cm}

Finally, let us complete the proof of Theorem 1. Let us write down
explicitly the divisor $\pi^*_+Y^+$ in terms of the partial
resolution $\sigma$. Let ${\cal E}^+$ be the set of all
exceptional divisors of the morphism $\sigma$, the image of which
on $V'$ is a divisor and covers the base $S'$. Therefore, ${\cal
E}^+$ can be identified with a subset of the set ${\cal E}$. In
the course of the proof of Proposition 1.3 we established that
$$
{\cal M}^+={\cal M}\cap{\cal E}^+\neq\emptyset.
$$
Now we write
$$
\pi^*_+Y^+=\pi^*Y-\sum_{E\in{\cal E}^+}\varepsilon_+(E)E+\Xi^+.
$$
Besides, we have
$$
K^+=\sigma^*K_V+\sum_{E\in{\cal E}^+}a_+(E)E+\Xi_K,
$$
where all coefficients $a_+(E)$ are positive and the divisor
$\Xi_K$ is effective, pulled back from the base $S^+$ and the
image of each of its irreducible component on $V'$ has codimension
at least 2, so that the general curve $C^+\in{\cal C}^+$ does not
intersect the support of the divisor $\Xi_K$. Let $C\in{\cal C}$
be its image on the original variety $V$ and
$\overline{C}=\pi(C)\in\overline{\cal C}$ the projection of the
curve $C$ on the base $S$. For a general divisor $D\in\Sigma$ and
its strict transform $D^+\in\Sigma^+$ on $V^+$ the
scheme-theoretic intersection
$(D^+\circ\pi^{-1}_+(\overline{C^+}))$ is well defined, it is an
effective cycle of dimension $\delta=\mathop{\rm dim}F$ on $V^+$.
For its numerical class we have the presentation
\begin{equation}\label{16.05.2014.1}
\begin{array}{c}
\displaystyle
(D^+\circ\pi^{-1}_+(\overline{C^+}))\sim-n(\sigma^*K_V\circ\pi^{-1}_+
(\overline{C^+}))+ \\
\\
\displaystyle +\left(\left[\sum_{E\in{\cal
E}^+}(-na_+(E)-\varepsilon_+(E))E\right]\cdot C^+\right)F.
\end{array}
\end{equation}
Since $(\overline{C}\cdot Y)\geqslant 0$ and
$(C^+\cdot\pi^*_+Y^+)< 0$, we have
$$
\left(-\left[\sum_{E\in{\cal E}^+}\varepsilon_+(E)E\right] \cdot
C^+\right)< 0,
$$
so that in the formula (\ref{16.05.2014.1}) the intersection of
the divisor in square brackets with $C^+$ is negative. Therefore,
$$
\sigma_*(D^+\circ\pi^{-1}_+(\overline{C^+}))\sim-n(K_V\circ\pi^{-1}
(\overline{C}))+bF,
$$
where $b<0$. Since on the left we have an effective cycle of
dimension $\delta$ on $V$, we obtain a contradiction with the
condition (iii) of our theorem. Proof of Theorem 1 is complete.
Q.E.D.


\section{Varieties of general position}

In this section we state the explicit local conditions of general
position for the double spaces (subsection 2.1) and hypersurfaces
(subsection 2.2), defining the sets ${\cal W}_{\rm reg}\subset
{\cal W}$ and ${\cal F}_{\rm reg}\subset {\cal F}$. In subsection
2.1 we prove Theorem 2. In subsection 2.3-2.5 we prove a part of
the claim of Theorem 4: the estimate for the codimension of the
complement ${\cal F}\setminus {\cal F}_{\rm reg}$; in subsection
2.5 we also consider some immediate geometric implications of the
conditions of general position.\vspace{0.3cm}

{\bf 2.1. The double spaces of general position.} The open subset
${\cal W}_{\rm reg}\subset {\cal W}$ of hypersurfaces of degree
$2M$ in ${\mathbb P}={\mathbb P}^M$ is defined by local
conditions, which a hypersurface $W\in {\cal W}_{\rm reg}$ must
satisfy at {\it every} point $o\in W$. These conditions depend on
whether the point $o\in W$ is non-singular or
singular.\vspace{0.1cm}

First, let us consider the condition of general position for a
{\bf non-singular point} $o\in W$. Let $(z_1,\dots,z_M)$ be a
system of affine coordinates with the origin at the point $o$ and
$$
w=q_1+q_2+\dots +q_{2M}
$$
the affine equation of the branch hypersurface $W$, where the
polynomials $q_i(z_*)$ are homogeneous of degree $i=1,\dots,2M$.
At a non-singular point $o\in W$ (that is, $q_1\not\equiv 0$) the
hypersurface $W$ must satisfy the condition\vspace{0.1cm}

(W1) the rank of the quadratic form $q_2|_{\{q_1=0\}}$ is at least
2.\vspace{0.1cm}

{\bf Proposition 2.1.} {\it Violation of the condition (W1)
imposes on the coefficients of the quadratic form $q_2$ (with the
linear form $q_1$ fixed)
$$
\frac{(M-2)(M-1)}{2}
$$
independent conditions.}\vspace{0.1cm}

{\bf Proof} is obvious. Q.E.D.\vspace{0.1cm}

Now let us consider the condition of general position for a {\bf
singular point} $o\in W$. Let
$$
w=q_2+q_3+\dots +q_{2M}
$$
be the affine equation of the branch hypersurface $W$ with respect
to a system of affine coordinates $(z_1,\dots,z_M)$ with the
origin at the point $o$. At a singular point $o$ the hypersurface
$W$ must satisfy the condition\vspace{0.1cm}

(W2) the rank of the quadratic form $q_2$ is at least
4.\vspace{0.1cm}

{\bf Proposition 2.2.} {\it Violation of the condition (W2)
imposes on the coefficients of the quadratic form $q_2$
$$
\frac{(M-2)(M-1)}{2}
$$
independent conditions.}\vspace{0.1cm}

{\bf Proof} is obvious Q.E.D.\vspace{0.1cm}

Now we define the subset ${\cal W}_{\rm reg}\subset {\cal W}$,
requiring that $W\in {\cal W}_{\rm reg}$ satisfies the condition
(W1) at every non-singular and the condition (W2) at every
singular point. Obviously, ${\cal W}_{\rm reg}\subset {\cal W}$ is
a Zariski open subset (possibly, empty).\vspace{0.1cm}

{\bf Proposition 2.3.} {\it The following estimate holds:}
$$
\mathop{\rm codim}(({\cal W}\setminus {\cal W}_{\rm reg}) \subset
{\cal W})\geqslant \frac{(M-4)(M-1)}{2}.
$$

{\bf Proof}  is obtained by the standard arguments, see
\cite[Chapter 3]{Pukh13a}: one considers the incidence subvariety
$$
{\cal I}=\{(o,W)\,|\, o\in W\}\subset {\mathbb P}\times{\cal W};
$$
for a fixed point $o\in {\mathbb P}$ the codimension of the set of
hypersurfaces ${\cal W}_{\rm non-reg}(o)$, containing that point
and non-regular in it, is given by Propositions 2.1 and 2.2 (in
the singular case $M$ more independent conditions are added as
$q_1\equiv 0$). After that one computes the dimension of the set
$$
{\cal I}_{\rm non-reg}=\mathop{\bigcup}\limits_{o\in{\mathbb P}}
\{o\}\times {\cal W}_{\rm non-reg}(o)
$$
and considers the projection onto ${\cal W}$. This completes the
proof. Q.E.D.\vspace{0.1cm}

Obviously, for any hypersurface $W\in {\cal W}_{\rm reg}$ the
double cover $F\to {\mathbb P}$, branched over $W$, is an
irreducible algebraic variety. Moreover, by the condition (W2) the
variety $F$ belongs to the class of varieties with quadratic
singularities of rank at least 5 \cite{EP}. Recall that a variety
${\cal X}$ is a variety with quadratic singularities of rank at
least $r$, if in a neighborhood of every point $o\in {\cal X}$ the
variety ${\cal X}$ can be realized as a hypersurface in a
non-singular variety ${\cal Y}$, and the local equation ${\cal X}$
at the point $o$ is of the form $\beta_1(u_*)+\beta_2(u_*)+\dots
=0$, where $(u_*)$ is a system of local parameters at the point
$o\in {\cal Y}$, and either $\beta_1\not\equiv 0$, or
$\beta_1\equiv 0$ and $\mathop{\rm rk} \beta_2\geqslant r$. It is
clear that $\mathop{\rm codim} (\mathop{\rm Sing}{\cal X}\subset
{\cal X})\geqslant r-1$, so that the variety $F$ is factorial
\cite{CL}.\vspace{0.1cm}

Furthermore, it is easy to show (see \cite{EP}), that the class of
quadratic singularities of rank at least $r$ is stable with
respect to blow ups in the following sense. Let $B\subset {\cal
X}$ be an irreducible subvariety. Then there exists an open set
${\cal U}\subset {\cal Y}$, such that ${\cal U}\cap B\neq
\emptyset$, ${\cal U}\cap B$ is a non-singular algebraic variety
and for its blow up
$$
\sigma_B\colon {\cal U}^+\to {\cal U}
$$
we have that $({\cal X}\cap {\cal U})^+\subset {\cal U}^+$ is a
variety of quadratic singularities of rank at least $r$. In order
to see this, note the following simple fact: if ${\cal Z}\ni o$ is
a non-singular divisor on ${\cal Y}$, where ${\cal Z}\neq {\cal
X}$ and the scheme-theoretic restriction ${\cal X}|_{\cal Z}$ has
at the point $o$ a quadratic singularity of rank $l$, then ${\cal
X}$ has at the point $o$ a quadratic singularity of rank at least
$l$. Now if $B\not\subset \mathop{\rm Sing}{\cal X}$, then the
claim about stability is obvious. Therefore, we may assume that
$B\subset \mathop{\rm Sing}{\cal X}$. The open set ${\cal
U}\subset {\cal Y}$ can be chosen in such a way that $B\cap {\cal
U}$ is a non-singular subvariety and the rank of quadratic points
$o\in B\cap {\cal U}$ is constant and equal to $l\geqslant r$. But
then in the exceptional divisor ${\cal E}=\sigma_B^{-1}(B\cap
{\cal U})$ the divisor $({\cal X}\cap {\cal U})^+\cap {\cal E}$ is
a fibration into quadrics of rank $l$, so that $({\cal X}\cap
{\cal U})^+\cap {\cal E}$ has at most quadratic singularities of
rank at least $l$. Therefore, $({\cal X}\cap {\cal U})^+\subset
{\cal U}^+$ has quadratic singularities of rank at least $r$ as
well, according to the remark above. For an explicit analytic
proof, see \cite{EP}.\vspace{0.1cm}

The stability with respect to blow ups implies that the
singularities of the variety $F$ are terminal (for the particular
case of one blow up it is obvious: the discrepancy of an
irreducible exceptional divisor $({\cal X}\cap {\cal U})^+\cap
{\cal E}$ with respect to ${\cal X}$ is positive; every
exceptional divisor over ${\cal X}$ can be realized by a sequence
of blow ups of the centres). Finally, $F$ satisfies the condition
(ii) of Theorem 1, that is, the condition of divisorial
canonicity, see the proof of part (ii) of Theorem 2 in
\cite{Pukh05} and Theorem 4 in \cite{Pukh09b}. This completes the
proof of Theorem 2. Q.E.D.\vspace{0.3cm}


{\bf 2.2. Fano hypersurfaces of general position.} As in the case
of double space, the open subset ${\cal F}_{\rm reg}\subset {\cal
F}$ of hypersurfaces of degree $M$ in ${\mathbb P}={\mathbb P}^M$
is defined by the local conditions, which a hypersurface $F\in
{\cal F}_{\rm reg}$ should satisfy at every point $o\in F$. Again
these conditions are different for non-singular and singular
points $o\in F$. Consider first the conditions of general position
for a {\bf non-singular point} $o\in F$.\vspace{0.1cm}

Let $(z_1,\dots,z_M)$ be a system of affine coordinates with the
origin at the point $o$ and
$$
w=q_1+q_2+q_3+\dots +q_{M}
$$
the affine equation of the hypersurface $F$, where the polynomials
$q_i(z_*)$ are homogene\-ous of degree $i=1,\dots,M$. Here is the
list of conditions of general position, which a hypersurface $F$
should satisfy at a non-singular point $o$.\vspace{0.1cm}

(R1.1) The sequence
$$
q_1,q_2,\dots ,q_{M-1}
$$
is regular in the local ring ${\cal O}_{o,{\mathbb P}}$, that is,
the system of equations
$$
q_1=q_2=\dots =q_{M-1}=0
$$
defines a one-dimensional subset, a finite set of lines in
${\mathbb P}$, passing through the point $o$. In particular,
$q_1\not\equiv 0$.\vspace{0.1cm}

The equation $q_1=0$ defines the tangent space $T_oF$ (which we,
depending on what we need, will consider either a linear subspace
in ${\mathbb C}^M$, or as its closure, a hyperplane in ${\mathbb
P}$). Now set ${\bar q}_i=q_i|_{\{q_1=0\}}$ for $i=2,\dots,M$:
these are polynomials on the linear space $T_o F\cong {\mathbb
C}^{M-1}$. The condition (R1.1) means the regularity of the
sequence
$$
{\bar q}_2, {\bar q}_3,\dots , {\bar q}_{M-1}.
$$
Such form is more convenient for estimating the codimension of the
set of hypersurfa\-ces which do not satisfy the regularity
condition.\vspace{0.1cm}

(R1.2) The quadratic form ${\bar q}_2$ on the space $T_o F$ is of
rank at least 6, and the linear span of every irreducible
component of the closed algebraic set
$$
\{q_1=q_2=q_3=0\}
$$
in ${\mathbb C}^M$ is the hyperplane $\{q_1=0\}$, that is, the
tangent hyperplane $T_o F$.\vspace{0.1cm}

An equivalent wording of this condition: every irreducible
component of the closed set $\{ {\bar q}_2={\bar q}_3=0 \}$ in
${\mathbb P}^{M-2}={\mathbb P}(\{q_1=0\})$ is
non-degenerate.\vspace{0.1cm}

(R1.3) For any hyperplane $P\subset {\mathbb P}$, $P\ni o$,
different from the tangent hyperplane $T_o F\subset {\mathbb P}$,
the algebraic cycle of scheme-theoretic intersection of
hyperplanes $P$, $T_o F$, the projective quadric
$\overline{\{q_2=0\}}\subset {\mathbb P}$ and $F$, that is, the
cycle,
$$
(P\circ \overline{\{q_1=0\}} \circ \overline{\{q_2=0\}} \circ F),
$$
is irreducible and reduced. (The line above means the closure in
${\mathbb P}$ and the operation $\circ$ of taking the cycle of
scheme-theoretic intersection is considered here on the space
${\mathbb P}$, too.)\vspace{0.1cm}

Now let us consider the conditions of general position for a {\bf
singular point} $o\in F$. Let $(z_1,\dots,z_M)$ be a system of
affine coordinates with the origin at the point $o$ and
$$
f=q_2+q_3+\dots+q_M
$$
the affine equation of the hypersurface $F$, where the polynomials
$q_i(z_*)$ are homogene\-ous of degree $i=2,\dots, M$. Let us list
the conditions of general position which must be satisfied for the
hypersurface $F$ at a singular point $o$.\vspace{0.1cm}

(R2.1) For any linear subspace $\Pi\subset{\mathbb C}^M$ of
codimension $c\in\{0,1,2\}$ the sequence
\begin{equation}\label{15.05.2014.2}
q_2|_{\Pi},\dots,q_{M-c}|_{\Pi}
\end{equation}
is regular in the ring ${\cal O}_{o,\Pi}$, that is, the system of
equations
$$
q_2|_{{\mathbb P}(\Pi)}=\dots=q_{M-c}|_{{\mathbb P}(\Pi)}=0
$$
defines in the space ${\mathbb P}(\Pi)\cong{\mathbb P}^{M-c-1}$ a
finite set of points.\vspace{0.1cm}

(R2.2) The quadratic form $q_2(z_*)$ is of rank at least
8.\vspace{0.1cm}

(R2.3) Now let us consider $(z_1,\dots,z_M)$ as homogeneous
coordinates $(z_1:\dots :z_M)$ on ${\mathbb P}^{M-1}$. The divisor
$$
\{q_3|_{\{q_2=0\}}=0\}
$$
on the quadric $\{q_2=0\}$ is not a sum of three (not necessarily
distinct) hyperplane sections of this quadric, taken from the same
linear pencil.\vspace{0.1cm}

Now arguing in the word for word the same way as in subsection ï.
2.1, we conclude that any hypersurface $F\in {\cal F}_{\rm reg}$
is an irreducible projective variety with factorial terminal
singularities. Obviously, $K_F=-H_F$ and $\mathop{\rm Pic} F =
{\mathbb Z} H_F$, where $H_F$ is the class of a hyperplane section
$F\subset {\mathbb P}$, that is, $F$ is a Fano variety of index
one. In order to prove Theorem 4, we have to show the following
two facts:\vspace{0.1cm}

--- the inequality (\ref{14.06.2014.1}),\vspace{0.1cm}

--- the divisorial log-canonicity of the hypersurface
$F\in {\cal F}_{\rm reg}$, that is, the condition (ii) of Theorem
1 for the variety $F$.\vspace{0.1cm}

These two tasks are dealt with in the remaining part of this
section and \S 3, respectively.\vspace{0.3cm}


{\bf 2.3. The conditions of general position at a non-singular
point.} Let $o\in F$ be a non-singular point. Fix an arbitrary
non-zero linear form $q_1$ and consider the affine space of
polynomials
$$
q_1+{\cal P}^{\rm sing}=\{q_1+q_2+\dots +q_M\},
$$
where ${\cal P}^{\rm sing}$ is the space of polynomials of the
form $f=q_2+q_3+\dots +q_M$. Let ${\cal P}_i\subset \{q_1+{\cal
P}^{\rm sing}\}$, $i=1,2,3$, be the closures of the subsets,
consisting of such polynomials $f$, which do not satisfy the
condition (R1.i), respectively. Set
$$
c_i=\mathop{\rm codim} ({\cal P}_i\subset \{q_1+{\cal P}^{\rm
sing}\}).
$$

{\bf Proposition 2.4.} {\it For $M\geqslant 8$ the following
equality holds:}
$$
\min\{c_1,c_2,c_3\}=c_2=\frac{(M-6)(M-5)}{2}.
$$

{\bf Proof} is easy to obtain by elementary methods. First of all,
by Lemma 2.1, shown below (where one must replace $M$ by $(M-1)$),
we obtain
$$
c_1=\frac{(M-1)(M-2)}{2} + 2.
$$
Furthermore, a violation of the condition $\mathop{\rm rk} {\bar
q}_2\geqslant 6$ imposes on the coefficients of the quadratic form
$q_2$
$$
\frac{(M-6)(M-5)}{2}< c_1
$$
independent conditions. Assuming the condition $\mathop{\rm rk}
{\bar q}_2\geqslant 6$ to be satisfied, we obtain that the quadric
$\{{\bar q}_2=0\}$ is factorial. It is easy to check that
reducibility or non-reducedness of the divisor $q_3|_{\{{\bar
q}_2=0\}}$ on this quadric gives
$$
\frac{M^3-6M^2-7M+54}{6}>\frac{(M-6)(M-5)}{2}
$$
independent conditions on the coefficients of the cubic form
$q_3$.\vspace{0.1cm}

Finally, let us consider a hyperplane $P\neq T_oF$ and the
quadratic hypersurface
$$
q_2|_{P\cap \{q_1=0\}}=0.
$$
Its rank is at least 5, so it is still factorial. Let us estimate
from below the number of independent conditions, which are imposed
on the coefficients of the polynomials $q_3,\dots, q_M$ if the
condition (R1.3) is violated. Define the values $v(\mu)$,
$\mu=0,1,2,3$, by the table \vspace{0.1cm}

\begin{center}
\begin{tabular}{|l|l|l|l|l|}
\hline $\mu$ & 0 & 1 & 2 & 3 \\
\hline $v(\mu)$ & 0 & 1 & $M$ & $\frac12 M(M+1) -1$ \\
\hline
\end{tabular}
\end{center}
\vspace{0.1cm}

\noindent  and set
$$
f(j,\mu)={j+M-1\choose M-1} - {j+M-3\choose M-1} -v(\mu)+ v(\max
(0, \mu-2)).
$$
Now, using the factoriality of the quadric, we obtain the estimate
$$
\begin{array}{c}
c_3\geqslant f(M,3)-(M-2) - \\ \\
- \max\left[ \max\limits_{M-1\geqslant j\geqslant 2}
(f(j,2)+f(M-j,1)), \max\limits_{M-1\geqslant j \geqslant 3}
(f(j,3)+f(M-j,0))\right].
\end{array}
$$
An elementary check shows that the minimum of the right hand side
is strictly higher than $c_2$ (and for $M\to \infty$ grows
exponentially). Q.E.D. for the proposition.\vspace{0.3cm}


{\bf 2.4. The conditions of general position at a singular point.}
Recall that ${\cal P}^{\rm sing}$ is the space of polynomials of
the form
$$
f=q_2+q_3+\dots+q_M
$$
in the variables $z_*=(z_1,\dots,z_M)$, where $q_i(z_*)$ are
homogeneous of degree $i$. Let ${\cal P}^{\rm sing}_{\rm
reg}\subset{\cal P}^{\rm sing}$ be the subset of polynomials
satisfying the conditions (R2.1-R2.3).\vspace{0.1cm}

{\bf Proposition 2.5.} {\it The following estimate holds:}
$$
\mathop{\rm codim}(\overline{({\cal{P}^{\rm sing}\backslash{\cal
P}^{\rm sing}_{\rm reg}})}\subset{\cal P}^{\rm
sing})=\frac{(M-7)(M-6)}{2}.
$$

{\bf Proof.} It is sufficient to show that violation of each of
the conditions (R2.1-R2.3) at the point $o=(0,\dots,0)$ separately
imposes on the polynomial $f$ at least $(M-7)(M-6)/2$ independent
conditions. It is easy to check that violation of the condition
(R2.2) imposes on the coefficients of the quadratic form
$q_2(z_*)$ precisely $(M-7)(M-6)/2$ independent conditions.
Therefore, considering the condition (R2.3), we may assume that
the condition (R2.2) is satisfied; in particular, the quadric
$\{q_2=0\}$) is factorial and violation of the condition (R2.3)
imposes on the coefficients of the cubic form $q_3(z_*)$ (with the
polynomial $q_2$ fixed)
$$
M\frac{M^2+3M-16}{6}\geqslant\frac{(M-7)(M-6)}{2}
$$
independent conditions for $M\geqslant 4$. It remains to consider
the case when the condition (R2.1) is violated.\vspace{0.1cm}

{\bf Lemma 2.1.} {\it Violation of the condition (R2.1) for one
value of the parameter $c=0$ imposes on the coefficients of the
polynomial $f$
\begin{equation}\label{15.05.2014.1}
\frac{M(M-1)}{2}+2
\end{equation}
independent conditions.}\vspace{0.1cm}

{\bf Proof} is obtained by the standard methods \cite[Chapter
3]{Pukh13a}. We just remind the scheme of arguments. Fix the first
moment when the sequence of polynomials $q_2,\dots,q_M$ becomes
non-regular: assume that the regularity is first violated for
$q_k$, that is, the closed set $\{q_2=\dots=q_{k-1}=0\}$ has the
``correct'' codimension $k-2$ and $q_k$ vanishes on one of the
components of that set. For $k\leqslant M-1$ we apply the method
of \cite{Pukh98b} and obtain that violation of the regularity
condition imposes on the coefficients of the polynomial $f$ at
least
$$
{M+1\choose k}\geqslant\frac{(M+1)M}{2}
$$
independent conditions; the right hand side of the last inequality
is strictly higher than (\ref{15.05.2014.1}), which is what we
need.\vspace{0.1cm}

Let us consider the last option:
$$
\{q_2=\dots=q_{M-1}=0\}\subset{\mathbb P}^{M-1}
$$
is a one-dimensional closed set and $q_M$ vanishes on one of its
irreducible components, say $B$. The case when $B\subset{\mathbb
P}^{M-1}$ is a line is a special one: it is easy to check that
vanishing on a line in ${\mathbb P}^{M-1}$ imposes on the
polynomials $q_2,\dots,q_M$ in total precisely
(\ref{15.05.2014.1}) independent conditions. Therefore, we may
assume that $B$ is not a line, that is, $\mathop{\rm dim} B
<\langle B\rangle =k\geqslant 2$. Now we apply the method
suggested in \cite{Pukh01}, fixing $k$ and the linear subspace
$\langle B\rangle$. To begin with, consider the case $k\leqslant
M-2$. In that case there are indices
$$
i_1,\dots,i_{k-1}\in\{2,\dots,M-1\},
$$
such that the restrictions $q_{i_1}|_{\langle B\rangle,}\dots,
q_{i_{k-1}}|_{\langle B\rangle}$ form a good sequence and $B$ is
one of its associated subvarieties (see \cite[Sec.3, Proposition
4]{Pukh01}, the details of this procedure are described in the
proof of the cited proposition). Taking into account that
$B\subset\langle B\rangle$ is by construction a non-degenerate
curve, we see that decomposable polynomials of the form $l_1\dots
l_a$, where $l_i$ are linear forms on $\langle
B\rangle\cong{\mathbb P}^k$, can not vanish on $B$. This gives
$jk+1$ independent conditions for each of the polynomials $q_j$
for $j\not\in\{i_1,\dots,i_{k-1}\}$, so that in total we get at
least
$$
\frac{k(M-k)(M-k+1)}{2}+M-2k-1
$$
independent conditions for these polynomials (the minimum is
attained for $i_1=M-k+1,\dots$, $i_{k-1}=M-1$). Taking into
account the condition $q_M|_B\equiv 0$ and the dimension of the
Grassmanian of $k$-dimensional subspaces in ${\mathbb P}^{M-1}$,
we obtain at least
$$
M^2-kM+k^2-M+k+1
$$
independent conditions for $f$. It is easy to check that the last
number is not smaller than (\ref{15.05.2014.1}).\vspace{0.1cm}

Finally, if $k=M-1$, that is, $B$ is a non-degenerate curve in
${\mathbb P}^{M-1}$, then the condition $q_M|_B\equiv 0$ gives at
least $M(M-1)+1$ independent conditions for $q_M$. Proof of Lemma
2.1 is complete. Q.E.D.\vspace{0.1cm}

Now let us complete the proof of Proposition 2.5.\vspace{0.1cm}

For a {\it fixed} linear subspace $\Pi\subset{\mathbb C}^M$ of
codimension $c\in\{0,1,2\}$ violation of regularity of the
sequence (\ref{15.05.2014.2}) imposes on the polynomial $f$ at
least $(M-c)(M-c-1)/2+2$ independent conditions. Subtracting the
dimension of the Grassmanian of subspaces of codimension $c$ in
${\mathbb C}^M$, we get the least value $(M-3)(M-6)/2$ for $c=2$.
This completes the proof of Proposition 2.5. Q.E.D.\vspace{0.3cm}


{\bf 2.5. Estimating the codimension of the complement to the set
${\cal F}_{\rm reg}$.} Recall that $F\in {\cal F}_{\rm reg}$ if
and only if at every non-singular point $o\in F$ the conditions
(R1.1-3) are satisfied, and at every singular point $o\in F$ the
conditions (R2.1-3) are satisfied. Propositions 2.4 and 2.5 imply
the following fact.\vspace{0.1cm}

{\bf Proposition 2.6.} {\it The following estimate holds:}
$$
\mathop{\rm codim}(({\cal F}\setminus {\cal F}_{\rm reg})\subset
{\cal F})\geqslant \frac{(M-7)(M-6)}{2}-5.
$$

{\bf Proof} is completely similar to the proof of Proposition 2.3
and follows from Propositions 2.4 and 2.5.\vspace{0.1cm}

Now let us consider some geometric facts which follow immediately
from the conditions of general position. These facts will be
needed in \S 3 to exclude log maximal singularities. In
\cite{Pukh05} it was shown that for any effective divisor $D\sim
nH$ on $F$ (where we write $H$ in stead of $H_F$ to simplify the
notations) the pair $(F,\frac{1}{n}D)$ is canonical at
non-singular points $o\in F$. This fact will be used without
special references. Now let $D_2=\{q_2|_F=0\}$ be the first
hypertangent divisor, so that we have $D^+_2\in|2H-3E|$. Recall
that $E\subset{\mathbb P}^{M-1}$ is an irreducible quadric of rank
at least 8. Obviously, the divisor $D_2\in|2H|$ satisfies the
equality
$$
\frac{\mathop{\rm mult}_o}{\mathop{\rm deg}}D_2=\frac{3}{M}.
$$
Here and below the symbol $\mathop{\rm
mult}\nolimits_o/\mathop{\rm deg}$ means the ratio of multiplicity
at the point $o$ to the degree.\vspace{0.1cm}

{\bf Lemma 2.2.} {\it Let $P\subset F$ be the section of the
hypersurface $F$ by an arbitrary linear subspace in ${\mathbb P}$
of codimension two, containing the point $o$. Then the restriction
$D_2|_P$ is an irreducible reduced divisor on the hypersurface}
$P\subset{\mathbb P}^{M-2}$.\vspace{0.1cm}

{\bf Proof.} The variety $P$ has at most quadratic singularities
of rank at least 6 and for that reason it is factorial. Therefore,
reducibility or non-reducedness of the divisor $D_2|_P$ means,
that the equality $D_2|_P=H_1+H_2$ holds, where $H_i$ are possibly
coinciding hyperplane sections of $P$. By the condition (R2.2) the
equalities $\mathop{\rm mult}_oH_i=2$ hold. However, $\mathop{\rm
mult}_oD_2|_P=6$. Therefore, $D_2|_P$ can not break into two
hyperplane sections. Q.E.D. for the lemma.\vspace{0.1cm}

{\bf Proposition 2.7.} {\it The pair $(F,\frac12 D_2)$ has no non
log canonical singularities, the centre of which on $F$ contains
the point} $o$: $LCS(F,\frac12 D_2)\not\ni o$.\vspace{0.1cm}

{\bf Proof.} Assume the converse. In any case
$$
\mathop{\rm codim} (LCS\left(F,\frac12 D_2\right)\subset
F)\geqslant 6,
$$
so that consider the section $P\subset F$ of the hypersurface $F$
by a generic linear subspace of dimension 5, containing the point
$o$. Then the pair $(P,\frac12 D_2|_P)$ has the point $o$ as an
isolated centre of a non log canonical singularity. Let
$\sigma_P\colon P^+\to P$ be the blow up of the non-degenerate
quadratic singularity $o\in P$ so that $E_P=E\cap P^+$ is a
non-singular exceptional quadric in ${\mathbb P}^4$. Since
$$
\frac12(D_2|_P)^+\sim H_P-\frac32 E_P\quad\mbox{and}\quad
a(E_P,P)=2>\frac32
$$
(where $H_P$ is the class of a hyperplane section of
$P\subset{\mathbb P}^5$), the pair $(P^+,\frac12(D_2|_P)^+)$ is
not log canonical. The union $LCS(P^+,\frac12(D_2|_P)^+)$ of all
centres of non log canonical singularities of that pair,
intersecting $E_P$, is a connected closed subset of the
exceptional quadric $E_P$, every irreducible component $S_P$ of
which satisfies the inequality $\mathop{\rm
mult}_{S_P}(D_2|_P)^+\geqslant 3$. Coming back to the original
pair $(F,\frac12 D_2)$, we see that for some irreducible
subvariety $S\subset E$ the inequality $\mathop{\rm
mult}_SD^+_2\geqslant 3$ holds, where $S\cap P^+=S_P$, so that
$\mathop{\rm codim}(S\subset E)\in\{1,2,3\}$.\vspace{0.1cm}

However, the case $\mathop{\rm codim}(S\subset E)=3$ is
impossible: by the connectedness principle this equality means
that $S_P$ is a point, and then $S\subset E$ is a linear subspace
of codimension 3, which is impossible if $\mathop{\rm
rk}q_2\geqslant 8$ (a 7-dimensional non-singular quadric does not
contain linear subspaces of codimension 3).\vspace{0.1cm}

Consider the case $\mathop{\rm codim}(S\subset E)=2$. Let
$\Pi\subset E$ be a general linear subspace of maximal dimension.
Then $D^+_2|_{\Pi}$ is a cubic hypersurface that has multiplicity
3 along an irreducible subvariety $S_{\Pi}=S\cap\Pi$ of
codimension 2. Therefore, $D^+_2|_{\Pi}$ is a sum of three (not
necessarily distinct) hyperplanes in $\Pi$, containing the linear
subspace $S_{\Pi}\subset\Pi$ of codimension 2, and for that reason
$D^+_2|_E$ is a sum of three (not necessarily distinct) hyperplane
sections from the same linear pencil as well, and $S$ is the
intersection of the quadric $E$ and a linear subspace of
codimension 2. However, this is impossible by the condition
(R2.3).\vspace{0.1cm}

Finally, if $\mathop{\rm codim}(S\subset E)=1$, then $D^+_2|_E=3S$
is a triple hyperplane section of the quadric $E$, which is
impossible by the condition (R2.3).\vspace{0.1cm}

This completes the proof of Proposition 2.7. Q.E.D.\vspace{0.1cm}

Here is one more fact that will be useful later.\vspace{0.1cm}

{\bf Proposition 2.8.} {\it For any hyperplane section $\Delta\ni
o$ of the hypersurface $F$ the pair $(F,\Delta)$ is log
canonical.}\vspace{0.1cm}

{\bf Proof.} This follows from a well known fact (see, for
instance, \cite{Co00,Pukh09b}): if $(p\in X)$ is a germ of a
non-degenerate quadratic three-dimensional singularity,
$\sigma\colon \widetilde{X}\to X$ its resolution with the
exceptional quadric $E_X\cong{\mathbb P}^1\times{\mathbb P}^1$ and
$D_X$ a germ of an effective divisor such that $o\in D_X$ and
$\widetilde{D}_X\sim -\beta E_X$, then the pair
$(X,\frac{1}{\beta}D_X)$ is log canonical at the point $o$. Q.E.D.


\section{Exclusion of maximal singularities}

In this section we complete the proof of Theorem 4. The symbol $F$
stands for a fixed hypersurface of degree $M$ in ${\mathbb P}$,
satisfying the regularity conditions: $F\in {\cal F}_{\rm reg}$.
As we mentioned in \S 2, in \cite{Pukh05} it was shown that the
pair $(F, \frac{1}{n} D)$ has no maximal singularities, the centre
of which is not contained in the closed set $\mathop{\rm Sing} F$,
for every effective divisor $D\sim nH$. In \cite{EP} it was shown
that for any mobile linear system $\Sigma\subset |nH|$ the pair
$(F, \frac{1}{n} D)$ is canonical for a general divisor $D\in
\Sigma$, that is, $\Sigma$ has no maximal singularities.
Therefore, in order to complete the proof of Theorem 4 it is
sufficient to show that for any effective divisor $D\sim nH$ the
pair $(F, \frac{1}{n} D)$ is log canonical, and we may assume only
those log maximal singularities, the centre of which is contained
in $\mathop{\rm Sing} F$.\vspace{0.1cm}

In subsection 3.1 we carry out preparatory work: by means of the
technique of hypertangent divisors we obtain estimates for the
ratio $\mathop{\rm mult}\nolimits_o/\mathop{\rm deg}$ for certain
classes of irreducible subvarieties of the hypersurface $F$. After
that we fix a pair $(F, \frac{1}{n} D)$ and assume that it is not
log canonical. The aim is to bring this assumption to a
contradiction. Let $B^*\subset \mathop{\rm Sing} F$ be the centre
of the log maximal singularity of the divisor $D$, $o\in B^*$ a
point of general position, $F^+\to F$ its blow up, $D^+$ the
strict transform of the divisor $D$. In subsection 3.2 we study
the properties of the pair $(F^+, \frac{1}{n} D^+)$: we show that
this pair has a non log canonical singularity, the centre of which
is a subvariety of the exceptional divisor of the blow up of the
point $o$. After that in subsections 3.2 and 3.3 we show that this
is impossible, which completes the proof of Theorem
4.\vspace{0.3cm}

{\bf 3.1. The method of hypertangent divisors.} Fix a singular
point $o\in F$,  a system of coordinates $(z_1,\dots,z_M)$ on
${\mathbb P}$ with the origin at that point and the equation
$f=q_2+\dots q_M$ of the hypersurface $F$.\vspace{0.1cm}

{\bf Proposition 3.1.} {\it Assume that the variety $F$ satisfies
the conditions (R2.1, R2.2) at the singular point $o$. Then the
following claims hold.\vspace{0.1cm}

{\rm (i)} For every irreducible subvariety of codimension 2
$Y\subset F$ the following inequality holds:
$$
\frac{\mathop{\rm mult}_o}{\mathop{\rm deg}}Y\leqslant\frac{4}{M}.
$$

{\rm (ii)} Let $\Delta\ni o$ be an arbitrary hyperplane section of
the hypersurface $F$. For every prime divisor $Y\subset \Delta$
the following inequality holds:
$$
\frac{\mathop{\rm mult}_o}{\mathop{\rm deg}}Y\leqslant\frac{3}{M}.
$$

{\rm (iii)} Let $P\ni o$ be the section of the hypersurface $F$ by
an arbitrary linear subspace of codimension two. For every prime
divisor $Y\subset P$ the following inequality holds:}
$$
\frac{\mathop{\rm mult}_o}{\mathop{\rm deg}}Y\leqslant\frac{4}{M}.
$$

{\bf Proof} is obtained by means of the method of hypertangent
divisors \cite[Chapter 3]{Pukh13a}. For $k=2,\dots,M-1$ let
$$
\Lambda_k=\left|\sum^k_{i=2}s_{k-i}(q_2+\dots+q_i)|_F=0\right|
$$
be the $k$-{\it th hypertangent linear system}, where $s_j(z_*)$
are all possible homogeneous polynomials of degree $j$. For the
blow up $\sigma\colon F^+\to F$ of the point $o$ with the
exceptional divisor $E=\sigma^{-1}(o)$, naturally realized as a
quadric in ${\mathbb P}^{M-1}$, we have
$$
\Lambda^+_k\subset |kH-(k+1)E|
$$
(where $\Lambda^+_k$ is the strict transform of the system
$\Lambda_k$ on $F^+$). Let $D_k\in\Lambda_k$, $k=2,\dots,M-1$ be
general hypertangent divisors.\vspace{0.1cm}

Let us show the claim (i). By the condition (R2.1) the equality
\begin{equation}\label{05.05.2014.1}
\mathop{\rm codim}\nolimits_o(\mathop{\rm Bs}\Lambda_k\subset
F)=k-1
\end{equation}
holds, where the symbol $\mathop{\rm codim}_o$ means the
codimension in a neighborhood of the point $o$; therefore,
$$
Y\cap D_4\cap D_5\cap\dots\cap D_{M-1}
$$
in a neighborhood of the point $o$ is a closed one-dimensional
set. We construct a sequence of irreducible subvarieties
$Y_i\subset F$ of codimension $i$: $Y_2=Y$ and $Y_{i+1}$ is an
irreducible component of the effective cycle $(Y_i\circ D_{i+2})$
with the maximal value of the ratio $\mathop{\rm
mult}_o/\mathop{\rm deg}$. The cycle $(Y_i\circ D_{i+2})$ every
time is well defined, because by the inequality
(\ref{05.05.2014.1}) we have $Y_i\not\subset D_{i+2}$ for a
general hypertangent divisor $D_{i+2}$. By the construction of
hypertangent linear system, at every step of our procedure the
inequality
$$
\frac{\mathop{\rm mult}_o}{\mathop{\rm
deg}}Y_{i+1}\geqslant\frac{i+3}{i+2}\cdot\frac{\mathop{\rm
mult}_o}{\mathop{\rm deg}}Y_i
$$
holds, so that for the curve $Y_{M-2}$ we have the estimates
$$
1\geqslant\frac{\mathop{\rm mult}_o}{\mathop{\rm
deg}}Y_{M-2}\geqslant \frac{\mathop{\rm mult}_o}{\mathop{\rm
deg}}Y\cdot\frac{5}{4}\cdot\frac{6}{5}\cdot\dots\cdot\frac{M}{M-1},
$$
which implies the claim (i).\vspace{0.1cm}

Let us prove the claim (ii). By Lemma 2.2, the divisor
$D_2|_{\Delta}$ is irreducible and reduced, and by the condition
(R2.1) it satisfies the equality
$$
\frac{\mathop{\rm mult}_o}{\mathop{\rm
deg}}D_2|_{\Delta}=\frac{3}{M}.
$$
Therefore, we may assume that $Y\neq D_2|_{\Delta}$, so that
$Y\not\subset D_2$ and the effective cycle of codimension two
$(Y\circ D_2)$ on $\Delta$ is well defined and satisfies the
inequality
$$
\frac{\mathop{\rm mult}_o}{\mathop{\rm deg}}(Y\circ
D_2)\geqslant\frac32\cdot\frac{\mathop{\rm mult}_o}{\mathop{\rm
deg}}Y.
$$
Let $Y_2$ be an irreducible component of that cycle with the
maximal value of the ratio $\mathop{\rm mult}_o/\mathop{\rm deg}$.
Applying to $Y_2$ the technique of hypertangent divisors in
precisely the same way as in the part (i) above, we see that by
the condition (R2.1) the intersection
$$
Y_2\cap D_4|_{\Delta}\cap D_5|_{\Delta}\cap\dots\cap
D_{M-2}|_{\Delta}
$$
in a neighborhood of the point $o$ is a one-dimensional closed
set, where $D_4\in\Lambda_4,\dots$, $D_{M-2}\in\Lambda_{M-2}$ are
general hypertangent divisors (note that the last hypertangent
divisor is $D_{M-2}$, and not $D_{M-1}$, as in the part (i),
because the dimension of $\Delta$ is one less than the dimension
of $F$ and the condition (R2.1) provides the regularity of the
truncated sequence $q_2|_{\Pi},\dots,q_{M-1}|_{\Pi}$, where in
this case $\Pi$ is a hyperplane, cutting out $\Delta$ on $F$).
Now, arguing in the word for word same way as in the proof of the
claim (i), we obtain the estimate
$$
1\geqslant\frac{\mathop{\rm mult}_o}{\mathop{\rm
deg}}Y\cdot\frac32\cdot\frac54\cdot\frac65\cdot\dots\cdot\frac{M-1}{M-2},
$$
which implies that
$$
\frac{\mathop{\rm mult}_o}{\mathop{\rm
deg}}Y\leqslant\frac{8}{3(M-1)}.
$$
For $M\geqslant 9$ the right hand side of the inequality does not
exceed $3/M$, which proves the claim (ii).\vspace{0.1cm}

Let us show the claim (iii). We argue in the word for word same
way as in the proof of the part (ii), with the only difference: in
order to estimate the multiplicity of the cycle $(Y\circ D_2|_P)$
at the point $o$ we use the hypertangent divisors
$$
D_4|_P,D_5|_P,\dots,D_{M-3}|_P
$$
(one less than above), so that we get the estimate
$$
1\geqslant\frac{\mathop{\rm mult}_o}{\mathop{\rm
deg}}Y\cdot\frac32\cdot\frac54\cdot\frac65\cdot\dots\cdot\frac{M-2}{M-3},
$$
which implies that
$$
\frac{\mathop{\rm mult}_o}{\mathop{\rm
deg}}Y\leqslant\frac{8}{3(M-2)}.
$$
For $M\geqslant 6$ the right hand side of the inequality does not
exceed $4/M$, which proves the claim (iii).\vspace{0.1cm}

Proof of Proposition 3.1 is complete. Q.E.D.\vspace{0.1cm}

Let us resume the proof of Theorem 4.\vspace{0.3cm}


{\bf 3.2. The blow up of a singular point.} Assume that the pair
$(F,\frac{1}{n}D)$ is not log canonical for some divisor
$D\in|nH|$, that is, for some prime divisor $E^*$ over $F$, that
is, a prime divisor $E^*\subset\widetilde{F}$, where
$\psi\colon\widetilde{F}\to F$ is some birational morphism,
$\widetilde{F}$ is non-singular and projective, the {\it log
Noether-Fano inequality} holds:
$$
\mathop{\rm ord}\nolimits_{E^*}\psi^* D> n(a(E^*)+1).
$$
By linearity of the inequality in $D$ and $n$ we may assume the
divisor $D$ to be prime. Let $B^*=\psi(E^*)\subset F$ be the
centre of the log maximal singularity $E^*$. We known that
$B^*\subset \mathop{\rm Sing} F$; in particular, $\mathop{\rm
codim}(B^*\subset F)\geqslant 7$. Let $o\in B^*$ be a point of
general position, $\varphi\colon F^+\to F$ its blow up, $E\subset
F^+$ the exceptional quadric.\vspace{0.1cm}

Consider the first hypertangent divisor $D_2\in|2H|$ at the point
$o$. By Lemma 2.2, the divisor $D_2$ is irreducible and reduced,
and by Proposition 2.7, the pair $(F,\frac12 D_2)$ is log
canonical at the point $o$. Therefore, $D\neq D_2$.\vspace{0.1cm}

{\bf Proposition 3.2.} {\it The following inequality holds}
$$
\mathop{\rm mult}\nolimits_oD\leqslant\frac{8}{3}n.
$$

{\bf Proof.} Consider the effective cycle $(D\circ D_2)$ of
codimension two. Obviously,
$$
\frac{\mathop{\rm mult}_o}{\mathop{\rm deg}}(D\circ
D_2)\geqslant\frac32\cdot \frac{\mathop{\rm mult}_o}{\mathop{\rm
deg}}D,
$$
however, by Proposition 3.1, (i), the left hand part of this
inequality does not exceed $(4/M)$. Since $\mathop{\rm deg}D=nM$,
Proposition 3.2 is shown. Q.E.D.\vspace{0.1cm}

Write down $D^+\sim nH-\nu E$, where $\nu\leqslant\frac43
n$.\vspace{0.1cm}

Let us consider the section $P$ of the hypersurface $F$ by a
general 5-dimensional linear subspace, containing the point $o$.
Let $P^+$ be the strict transform of $P$ on $F^+$ and $E_P=P^+\cap
E$ a non-singular three-dimensional quadric. Set also $D_P=D|_P$.
Obviously, the pair $(P,\frac{1}{n}D_P)$ has the point $o$ as an
isolated centre of a non log canonical singularity. Since
$a(E_P)=2$ and $D^+_P\sim nH_P-\nu E_P$ (where $H_P$ is the class
of a hyperplane section of the variety $P$), where
$\nu\leqslant\frac43n<2n$, the pair $(P^+,\frac{1}{n}D^+_P)$ is
not log canonical and the union $LCS_E(P^+,\frac{1}{n}D^+_P)$ of
centres of all non log canonical singularities of that pair,
intersecting $E_P$, is a connected closed subset of the
exceptional quadric $E_P$. Let $S_P$ be an irreducible component
of that set. Obviously, the inequality
$$
\mathop{\rm mult}\nolimits_{S_P}D^+_P> n
$$
holds. Furthermore, $\mathop{\rm codim}(S_P\subset
E_P)\in\{1,2,3\}$. Returning to the original pair
$(F,\frac{1}{n}D)$, we see that there is a non log canonical
singularity of the pair $(F^+,\frac{1}{n}D^+)$, the centre of
which is a subvariety $S\subset E$, such that $S\cap E_P=S_P$ and,
in particular, $\mathop{\rm codim}(S\subset
E)\in\{1,2,3\}$.\vspace{0.1cm}

Note at once that the case $\mathop{\rm codim}(S\subset E)=3$ is
impossible: by the connectedness principle in that case $S_P$ is a
point and for that reason $S$ is a linear subspace of codimension
3 on the quadric $E$ of rank at least 8, which is
impossible.\vspace{0.1cm}

It is not hard to exclude the case $\mathop{\rm codim}(S\subset
E)=1$, either. Assume that it does take place. Then the divisor
$S$ is cut out on $E$ by a hypersurface of degree $d_S\geqslant
1$. Let $H_E$ be the class of a hyperplane section of the quadric
$E$. The divisor $D^+|_E\sim\nu H_E$, so that
$$
\frac43n\geqslant\nu>nd_S
$$
and for that reason $S$ is a hyperplane section of the quadric
$E$. Let $\Delta\in|H|$ be the uniquely determined hyperplane
section of the hypersurface $F$, such that $\Delta\ni o$ and
$\Delta^+\cap E=S$. The pair $(F^+,\Delta^+)$ is log canonical and
for that reason $D\neq\Delta$. For the effective cycle
$(D\circ\Delta)$ of codimension two on $F$ we have
$$
\mathop{\rm mult}\nolimits_o(D\circ\Delta)\geqslant
2\nu+2\mathop{\rm mult}\nolimits_SD^+>4n,
$$
so that
$$
\frac{\mathop{\rm mult}_o}{\mathop{\rm
deg}}(D\circ\Delta)>\frac{4}{M},
$$
which contradicts Proposition 3.1. This excludes the case of a
divisorial centre.\vspace{0.3cm}


{\bf 3.3. The case of codimension two.} Starting from this moment,
assume that $\mathop{\rm codim}(S\subset E)=2$.\vspace{0.1cm}

{\bf Lemma 3.1.} {\it The subvariety $S$ is contained in some
hyperplane section of the quadric $E$.}\vspace{0.1cm}

{\bf Proof.} Since $\mathop{\rm mult}_SD^+>n$ and $D^+|_E\sim \nu
H_E$ with $\nu\leqslant\frac43n$, for every secant line $L\subset
E$ of the subvariety $S$ we have $L\subset D^+$. Let $\Pi\subset
E$ be a linear space of maximal dimension and of general position
and $S_{\Pi}=S\cap\Pi$. The secant lines of the closed set
$S_{\Pi}\subset\Pi$ of codimension two can not sweep out $\Pi$,
since $E\not\subset D^+$. Therefore, there are two options (see
\cite[Lemma 2.3]{Pukh09b}):\vspace{0.1cm}

1) the secant lines of the set $S_{\Pi}$ sweep out a hyperplane in
$\Pi$,\vspace{0.1cm}

2) $S_{\Pi}\subset\Pi$ is a linear subspace of codimension
two.\vspace{0.1cm}

In the first case the secant lines $L\subset E$ of the set $S$
sweep out a divisor on $E$, which ca only be a hyerplane section
of the quadric $E$. In the second case $S$ contains all its secant
lines and is a section of $E$ by a linear subspace of codimension
two. Q.E.D. for the lemma.\vspace{0.1cm}

As we have just shown, one of the two options takes place: either
there is a unique hyperplane section $\Lambda$ of the quadric $E$,
containing $S$ (Case 1), or $S=E\cap \Theta$, where $\Theta$ is a
linear subspace of codimension two (Case 2). Let us study them
separately.\vspace{0.1cm}

Assume that {\bf Case 1} takes place. Then  $S$ is cut out on
$\Lambda$ by a hypersurface of degree $d_S\geqslant2$. Set
$$
\mu=\mathop{\rm
mult}\nolimits_SD^+\quad\mbox{and}\quad\gamma=\mathop{\rm
mult}\nolimits_{\Lambda}D^+,
$$
where $\mu>n$ and $\mu\leqslant
2\nu\leqslant\frac83n$.\vspace{0.1cm}

{\bf Lemma 3.2.} {\it The following inequality holds:}
$$
\gamma\geqslant\frac{2\mu-\nu}{3}.
$$

{\bf Proof}  is easy to obtain in the same way as the short proof
of Lemma 3.5 in \cite[subsection 3.7]{Pukh00d}. Let
$L\subset\Lambda$ be a general secant line of the set $S$.
Consider the section $P$ of the hypersurface $F$ by a general
4-plane in ${\mathbb P}$, such that $P\ni o$ and $P^+\cap E$
contains the line $L$. Obviously, $o\in P$ is a non-degenerate
quadratic point and $E_P=P^+\cap E\cong{\mathbb P}^1\times{\mathbb
P}^1$ is a non-singular quadric in ${\mathbb P}^3$. Set
$D_P=D|_P$. Obviously, $\gamma=\mathop{\rm
mult}_LD^+_P$.\vspace{0.1cm}

Let $\sigma_L\colon P_L\to P^+$ be the blow up of the line $L$,
$E_L=\sigma^{-1}_L(L)$ the exceptional divisor; since ${\cal
N}_{L/P^+}\cong{\cal O}\oplus{\cal O}_L(-1)$, the exceptional
surface $E_L$ is a ruled surface of the type ${\mathbb F}_1$, so
that
$$
\mathop{\rm Pic}E_L={\mathbb Z}s\oplus{\mathbb Z}f,
$$
where $s$ and $f$ are the classes of the exceptional section and
the fibre, respectively. Furthermore, $E_L|_{E_L}=-s-f$. Let $D_L$
be the strict transform of $D^+_P$ on $P_L$. Obviously,
$$
D_L\sim n H_P-\nu E_P-\gamma E_L
$$
(where $H_P$ is the class of a hyperplane section of $P$), so that
$$
D_L|_{E_L}\sim\gamma s+(\gamma+\nu)f.
$$
On the other hand, $L$ is a general secant line of the set $S$ and
for that reason $L$ contains at least two distinct points $p,q\in
S$. Therefore, the divisor $D^+_P$ has at the points $p,q\in L$
the multiplicity $\mu$ and for that reason the effective 1-cycle
$D_L|_{E_L}$ contains the corresponding fibres $\sigma^{-1}_L(p)$
and $\sigma^{-1}_L(q)$ over those points with multiplicity
$(\mu-\gamma)$. Therefore,
$$
\gamma+\nu\geqslant 2(\mu-\gamma),
$$
whence follows the claim of the lemma. Q.E.D.\vspace{0.1cm}

Now let us consider the uniquely determined hyperplane section
$\Delta$ of the hypersur\-face $F\subset{\mathbb P}$, such that
$\Delta\ni o$ and $\Delta^+\cap E=\Lambda$. Set
$D_{\Delta}=D|_{\Delta}$. Write down
$$
D^+|_{\Delta^+}=D^+_{\Delta}+a\Lambda.
$$
Obviously,
$$
\mathop{\rm mult}\nolimits_oD_{\Delta}=2(\nu+a)\geqslant
2\nu+2\frac{2\mu-\nu}{3}=\frac43(\mu+\nu)>\frac83n.
$$
Since as we noted above, the subvariety $S$ is cut out on the
quadric $\Lambda$ by a hypersurfa\-ce of degree $d_S\geqslant 2$,
the divisor $D^+_{\Delta}\sim nH_{\Delta}-(\nu+a)\Lambda$ can not
contain $S$ with multiplicity higher than
$$
\frac{1}{d_S}(\nu+a)\leqslant\frac{\nu+a}{2}.
$$
Since the pair $(F^+,\frac{1}{n}D^+)$ has a non log canonical
singularity with the centre at $S$, the inversion of adjunction
implies that the pair
$\square=(\Delta^+,\frac{1}{n}(D^+_{\Delta}+a\Lambda))$ has a non
log canonical singularity with the centre at $S$ as well. Recall
that $\mathop{\rm codim}(S\subset\Delta^+)=2$. Consider the blow
up $\sigma_S\colon\widetilde{\Delta}\to\Delta^+$ of the subvariety
$S$ and denote by the symbol $E_S$ the exceptional divisor
$\sigma^{-1}_S(S)$. The following fact is well
known.\vspace{0.1cm}

{\bf Proposition 3.3.} {\it For some irreducible divisor
$S_1\subset E_S$, such that the projection $\sigma_S|_{S_1}$ is
birational, the inequality
\begin{equation}\label{03.05.2014.1}
\mathop{\rm mult}\nolimits_S(D^+_{\Delta}+a\Lambda)+\mathop{\rm
mult}\nolimits_{S_1}(\widetilde{D}_{\Delta}+a\widetilde{\Lambda})>2n
\end{equation}
holds, where  $\widetilde{D}_{\Delta}$ and  $\widetilde{\Lambda}$
are the strict transforms of $D^+_{\Delta}$ and $\Lambda$ on
$\widetilde{\Delta}$, respectively.}\vspace{0.1cm}

{\bf Proof:} see Proposition 9 in \cite{Pukh05}.\vspace{0.1cm}

Set $\mu_S=\mathop{\rm mult}_SD^+_{\Delta}$ and $\beta=\mathop{\rm
mult}_{S_1}\widetilde{D}_{\Delta}$. Consider first the case of
general position: $S_1\neq E_S\cap\widetilde{\Lambda}$. In that
case $S_1\not\subset\widetilde{\Lambda}$ and the inequality
(\ref{03.05.2014.1}) takes the following form:
$$
\mu_S+\beta+a>2n.
$$
Since $\mu_S\geqslant\beta$, the more so $2\mu_S+a>2n$. On the
other hand, we noted above that $2\mu_S\leqslant\nu+a$. As a
result, we obtain the estimate
$$
\nu+2a>2n.
$$
Therefore, $\mathop{\rm mult}_oD_{\Delta}>\nu+2n>3n$. However,
$D_{\Delta}\sim nH_{\Delta}$ is an effective divisor on the
hyperplane section $\Delta$ and by Proposition 3.1, (ii), it
satisfies the inequality
$$
\frac{\mathop{\rm mult}_o}{\mathop{\rm
deg}}D_{\Delta}\leqslant\frac{3}{M}.
$$
This contradiction excludes the case of general position.
Therefore, we are left with the only option:
$S_1=E_S\cap\widetilde{\Lambda}$.\vspace{0.1cm}

In that case the inequality (\ref{03.05.2014.1}) takes the
following form:
$$
\mu_S+\beta+2a>2n.
$$
This inequality is weaker than the corresponding estimate in the
case of general position, but as a compensation we obtain the
additional inequality
$$
2\mu_S+2\beta\leqslant\nu+a
$$
(the restriction $D^+_{\Delta}|_{\Lambda}$ is cut out by a
hypersurface of degree $(\nu+a)$ and contains the divisor $S\sim
d_SH_{\Lambda}$ with multiplicity at least $\mu_S+\beta$).
Combining the last two estimates, we obtain the inequality
$$
\nu+5a>4n,
$$
which implies that $5(\nu+a)>8n$ and so $\mathop{\rm
mult}_oD_{\Lambda}>\frac{16}{5}n$; as we mentioned above, this
contradicts Proposition 3.1, (ii). This completes the exclusion of
Case 1.\vspace{0.1cm}

Therefore, {\bf Case 2} takes place: $S=E\cap\Theta$, where
$\Theta\subset{\mathbb P}^{M-1}$ is a linear subspace of
codimension two. Let $P\subset F$ be the section of the
hypersurface $F$ by the linear subspace of codimension two in
${\mathbb P}$, which is uniquely determined by the conditions
$P\ni o$ and $P^+\cap E=S$. Furthermore, let $|H-P|$ be the pencil
of hyperplane sections of $F$, containing $P$. For a general
hyperplane section $\Delta\in|H-P|$ we have:\vspace{0.1cm}

\begin{itemize}

\item the divisor $D$ does not contain $\Delta$ as a component,
so that the effective cycle $(D\circ\Delta)=D_{\Delta}$ of
codimension two on $F$ is well defined; this cycle can be looked
at as an effective divisor $D_{\Delta}\in|nH_{\Delta}|$ on the
hypersurface $\Delta\subset{\mathbb P}^{M-1}$,

\item for the strict transform $D^+_{\Delta}$ on $F^+$ the
equality $\mathop{\rm mult}_SD^+_{\Delta}=\mathop{\rm mult}_SD^+$
holds.

\end{itemize}

Of course, the divisor $D_{\Delta}$ may contain $P$ as a
component. Write down $D_{\Delta}=G+aP$, where $a\in{\mathbb Z}_+$
and $G$ is an effective divisor that does not contain $P$ as a
component, $G\in|(n-a)H_{\Delta}|$. Obviously,
$G^+\sim(n-a)H_{\Delta}-(\nu-a)E_{\Delta}$, where
$E_{\Delta}=\Delta^+\cap E$, and, besides,
$$
\mathop{\rm mult}\nolimits_SG^+=\mathop{\rm
mult}\nolimits_SD^+-a>n-a.
$$
Set $m=n-a$. The effective cycle of codimension two $G_P=(F\circ
P)$ on $\Delta$ is well defined and can be considered as an
effective divisor $G_P\in|mH_P|$ on the hypersurface
$P\subset{\mathbb P}^{M-2}$. The divisor $G_P$ satisfies the
inequality
$$
\mathop{\rm mult}\nolimits_oG_P\geqslant 2(\nu-a)+ 2\mathop{\rm
mult}\nolimits_SG^+>4m.
$$
This is impossible by Proposition 3.1, (iii).\vspace{0.1cm}

Therefore, the assumption that the pair $(F,\frac{1}{n} D)$ is not
log canonical for some divisor $D\sim nH$, leads to a
contradiction.\vspace{0.1cm}

Proof of Theorem 4 is complete.

\begin{flushleft}
Department of Mathematical Sciences,\\
The University of Liverpool
\end{flushleft}

\noindent{\it pukh@liv.ac.uk}

\end{document}